\documentclass{ifacconf}

\usepackage{graphicx}      
\usepackage{natbib}        
\usepackage{mathtools}
\usepackage{tikz}
\usepackage{amsmath,amssymb,amsfonts}

\newcounter{MYtempeqncnt}
\begin{document}
\begin{frontmatter}

\title{Gaussian inference for data-driven state-feedback design of nonlinear systems\thanksref{footnoteinfo}} 

\thanks[footnoteinfo]{Frank Allgöwer thanks the funding by the Deutsche Forschungsgemeinschaft (DFG, German Research Foundation) under Germany's Excellence Strategy - EXC 2075 - 390740016 and under grant 468094890. The work of Thomas B. Schön is funded by \emph{Kjell och Märta Beijer Foundation} and by the project \emph{NewLEADS - New Directions in Learning Dynamical Systems} (contract number: 621-2016-06079), funded by the Swedish Research Council. Tim Martin and Frank Allgöwer also acknowledge the support by the Stuttgart Center for Simulation Science (SimTech).\\
\copyright2023 the authors. This work has been accepted to IFAC for publication under a Creative Commons Licence CC-BY-NC-ND.}

\author[First]{Tim Martin} 
\author[Second]{Thomas B. Schön} 
\author[First]{Frank Allgöwer}

\address[First]{University of Stuttgart, Institute for Systems Theory and Automatic Control, Germany, \{tim.martin, frank.allgower\}@ist.uni-stuttgart.de.}
\address[Second]{Uppsala University, Departement of Information Technology, Sweden, thomas.schon@it.uu.se.}


\begin{abstract}                
	Data-driven control of nonlinear systems with rigorous guarantees is a challenging problem as it usually calls for nonconvex optimization and often requires the knowledge of the true basis functions of the unknown system dynamics. To tackle these drawbacks, this work is based on a data-driven polynomial representation of general nonlinear systems exploiting Taylor polynomials. Thereby, we design state-feedback laws that render a known equilibrium globally asymptotically stable while operating with respect to a desired quadratic performance criterion. The calculation of the polynomial state feedback boils down to a sum-of-squares optimization problem, and hence to computationally tractable linear matrix inequalities. Moreover, we examine state-input data in presence of Gaussian noise by Bayesian inference to overcome the conservatism of deterministic noise characterizations from recent data-driven control approaches for Gaussian noise. 
\end{abstract}


\end{frontmatter}

\section{Introduction}

Controller design techniques \citep{Khalil} typically require a precise model of the system. However, applying first principles for modelling a system can be expensive in time and usually requires expert knowledge. To this end, interest in data-driven methods has risen, where a controller is received from measured trajectories. For example, system identification \citep{SysId} establishes an indirect procedure by first identifying a model from measurements and then applying model-based controller design tools. Here, closed-loop stability can only be guaranteed if the approximation error of the model is known which is however even for linear time-invariant systems an active research field \citep{Oymak}. Moreover, the amount of data required for identifying the dynamics can be larger than for stabilizing the system \citep{Informativity}.

Recent research includes direct data-driven approaches without first identifying an explicit model. For linear time-invariant (LTI) systems, \cite{PersisLinear} relies on the behavioral systems theory, \cite{vanWaarde} introduces a matrix S-lemma, and \cite{Berberich} uses a linear fraction representation to combine data and prior knowledge. As a step towards nonlinear systems, extensions for certain system classes as polynomial \citep{DePersis1} and rational systems \citep{Strasser} are examined. Data-driven approaches for general nonlinear systems include adaptive control \citep{Adaptive}, Koopman linearization \citep{Koopman}, feedback linearization \citep{FeedbackLin}, set-membership \citep{Milanese}, linearly parametrized models with known basis functions \citep{Sznaier} and \citep{DePersisNonlinearity}, and combining Gaussian processes and robust control techniques \citep{Umlauft} and \citep{Fiedler}. 

The mentioned methods mostly require the true basis functions to be known, lack on rigorous stability and performance guarantees, or require nonconvex optimization. To this end, we establish in this work a state-feedback design by the data-based representation of general nonlinear systems using Taylor polynomials (TP) from \cite{MartinTP2}. Thereby, a single sum-of-squares (SOS) synthesis condition is determined and thus leads to computationally appealing linear matrix inequalities (LMI). Since we first determine from data a suitable representation of the nonlinear system dynamics and its uncertainty to design a robust controller in a second step, the presented data-driven controller design is indirect similar to Koopman, set-membership, and Gaussian process approaches. Note that \cite{MartinTP2} tackles the problem of verifying dissipativity properties, which is structurally easier to solve than a controller design. Moreover, we consider here the TP representation subject to Gaussian noise instead of a deterministic noise characterization. Therefore, the presented work is in line with \cite{Umlauft}, \cite{Fiedler}, and \cite{Umenberger} where uncertainty inferences from probabilistic machine learning techniques are utilized for a robust controller design. Furthermore, Gaussian noise is interesting if only an inaccurate deterministic bound on the noise is available, and hence leads to impractical inferences. Gaussian noise is also a common assumption in system identification \citep{SysId} such that recent data-driven results could be compared to system identification techniques in future work. 

We make several contributions in this paper. By the extension of \cite{MartinTP2} to Gaussian noise, we generalize the data-based representation for LTI systems from \cite{Umenberger} to general nonlinear systems. At the same time, we consider not only the Bayesian treatment as in \cite{Umenberger} but also the so-called Frequentist treatment \citep{Bishop} which can be directly connected to the results for deterministic noise characterizations. Moreover, we show how prior knowledge of the system dynamics for the Bayesian inference from \cite{Umenberger} can be exploited to improve its accuracy. In particular, this plays a crucial role when applying TP representations for real data as indicated by \cite{MartinTP1}.

Further contributions are that we build our controller synthesis on the basis of the flexible LMI-based robust control framework of \cite{SchererLMI} to achieve a single SOS condition to determine a state feedback that guarantees to render a known equilibrium globally asymptotically stable. In contrast, the recent investigation in \cite{DePersis2} of a TP representation for designing state-feedback laws by Petersen's lemma only achieves asymptotic stability which additionally calls for an iterative approximation of the region of attraction. Furthermore, we allow for a synthesis with performance criteria, for instance, to reduce the required control input. Similar to \cite{Berberich}, we can also make use of prior knowledge on the system dynamics in the controller synthesis to reduce the number of required data and improve the control performance, which is essential for the application of the TP representation in practice \citep{MartinTP1}. Since a global representation of a nonlinear system by means of a single TP might have a large uncertainty inherent, we also provide a controller synthesis with local stability and performance guarantees. 

The paper is organized as follows. After providing some notation in Section~\ref{SecNotation}, we introduce our setup in Section~\ref{SecProbSetup}. In Section~\ref{SecPolySector}, the TP representation of nonlinear systems from \cite{MartinTP2} is recapped and extended to incorporate prior knowledge on the dynamics. Section~\ref{SecBayesianSetMem} presents two possibilities for a Gaussian inference on the unknown TP. Subsequent, the controller synthesis is considered in Section~\ref{SecConSyn}. Section~\ref{SecNumEx} compares both Gaussian inference schemes for the stabilization of an inverted pendulum in a numerical example.

\section{Notation}\label{SecNotation}

We denote the Euclidean norm of a vector $v\in\mathbb{R}^n$ by $||v||_2$ and the identity and the zero matrix of suitable dimensions by $I$ and $0$, respectively. The binomial coefficient $\begin{pmatrix} m\\n\end{pmatrix}$ is denoted by $C_{m,n}$, the Minkowski sum of two sets by $\oplus$, and the Kronecker product of two matrices by $\otimes$. For two matrices $A_1$ and $A_2$ of suitable dimensions, consider the abbreviation of a quadratic form $A_1^TA_2A_1=\star^T A_2\cdot A_1$ and 
\begin{equation*}
\begin{bmatrix}\begin{array}{c|c}
A_1 & 0 \\\hline 0 & A_2 
\end{array}	\end{bmatrix}=\text{diag}(A_1\big| A_2).
\end{equation*}
Furthermore, if a random vector $X$ is Gaussian distributed with mean $\mu$ and covariance matrix $\Xi\succ0$, then $X\sim\mathcal{N}(\mu,\Xi)$. 
$Q_k$ denotes the quantile function of the Chi-squared distribution with $k$ degrees of freedom, i.e., for a Chi-squared distributed random variable $Y$ with $k$ degrees of freedom, $p(Y\leq Q_k(\delta))=\delta$.

$\mathbb{R}[x]$ corresponds to the set of all real polynomials in $x=\begin{bmatrix}x_1 & \cdots &x_n\end{bmatrix}^T\in\mathbb{R}^n$
\begin{equation*}
p(x)=\sum_{\alpha\in\mathbb{N}^n,|\alpha|\leq d} a_\alpha x^\alpha,
\end{equation*}
with vectorial indices $\alpha^T=\begin{bmatrix}\alpha_1 & \cdots & \alpha_n\end{bmatrix}^T\in\mathbb{N}^n$, $|\alpha|=\alpha_1+\cdots+\alpha_n$, real coefficients $a_\alpha\in\mathbb{R}$, and monomials $x^\alpha=x_1^{\alpha_1}\cdots x_n^{\alpha_n}$. $d$ is called the degree of polynomial $p$. Furthermore, let $\mathbb{R}[x]^m$ denote the set of all $m$-dimensional polynomial vectors and $\mathbb{R}[x]^{m\times n}$ all $m\times n$ polynomial matrices which entries are polynomials from $\mathbb{R}[x]$. 
For a matrix $P\in\mathbb{R}[x]^{n\times n}$, if there exists a matrix $Q\in\mathbb{R}[x]^{{m}\times {n}}$ such that $P=Q^TQ$, then $P\in\text{SOS}[x]^n$ where $\text{SOS}[x]^n$ denotes the set of all ${n}\times {n}$ SOS matrices in $x$. Analogously for $n=1$, $P$ is called an SOS polynomial and $\text{SOS}[x]$ is the set of all SOS polynomials.

\section{Problem formulation}\label{SecProbSetup}

Throughout the paper, we study the continuous-time input-affine system
\begin{align}\label{TrueSystem}
\dot{x}(t)=f(x(t))+B(x(t))u(t)
\end{align}
with an unknown $k+1$ times continuously differentiable nonlinear function $f(x)=\begin{bmatrix}f_1(x)&\cdots&f_{n_x}(x)\end{bmatrix}^T:\mathbb{R}^{n_x}\rightarrow\mathbb{R}^{n_x}$ and unknown polynomial input matrix $B\in\mathbb{R}[x]^{n_x\times n_u}$. Without loss of generality, we assume that $f(0)=0$. Then, the goal of the paper is to derive a state-feedback law $u(x)$ that globally asymptotically stabilizes the known equilibrium $x=0$ of the unknown nonlinear system \eqref{TrueSystem} using the available noisy measurements
\begin{equation}\label{DataSystem}
\{(\dot{\tilde{x}}_i,\tilde{x}_i,\tilde{u}_i)_{i=1,\dots,S}\}
\end{equation}
with $\dot{\tilde{x}}_i=f(\tilde{x}_i)+B(\tilde{x}_i)\tilde{u}_i+\tilde{d}_i,i=1,\dots,S$. The unknown disturbance $\tilde{d}_i$ can take process noise and inaccurate estimates of $\dot{\tilde{x}}_i$ into account.

To achieve rigorous guarantees for the state feedback for a general nonlinear drift, further insights into $f$ are necessary. Indeed, inferring the dynamics \eqref{TrueSystem} at an unseen state $x$ is impossible from only a finite set of samples \eqref{DataSystem}. Thus, the following assumptions are appropriate.
\begin{assum}[\cite{MartinTP2}]\label{AssBoundDeri}
\hfill Upper\\ bounds $M_{i,\alpha}\geq0,i=1,\dots,n_x,|\alpha|=k+1,\alpha\in\mathbb{N}^{n_x},$ on the magnitude of each $(k+1)$-th order partial derivative are known, i.e.,
\begin{align*}
\bigg|\bigg|\frac{\partial^{k+1}f_i(x)}{\partial x^{\alpha}}\bigg|\bigg|_2\leq M_{i,\alpha},\quad \forall x\in\mathbb{R}^{n_x}. 
\end{align*}
\end{assum}
\begin{assum}[\cite{MartinPoly}]\label{AssBoundDegree}
An upper\\ bound on the degree of the polynomial matrix $B(x)$ is known.
\end{assum}
\begin{assum}\label{AssNoiseBound}
The disturbances $\tilde{d}_i,i=1,\dots,S,$ are independent and Gaussian $\tilde{d}_{i}\sim\mathcal{N}(0,\sigma^2I)$ with known standard deviation $\sigma$.
\end{assum} 

Since the information about the rate of variation of $f(x)$ according to Assumption~\ref{AssBoundDeri} is typically not available, \cite{MartinTP1} proposes a validation procedure to obtain reasonable bounds from noisy data, which were already applied in an experimental example. Moreover, note that the knowledge of Assumption~\ref{AssBoundDeri} for all $x\in\mathbb{R}^{n_x}$ might be restrictive. Hence, we also consider a local controller synthesis in Section~\ref{SecExtensions}. Assumption~\ref{AssNoiseBound} supposes Gaussian noise as common in system identification and which differ from the deterministic noise characterization in most data-driven robust controller results, e.g., \cite{vanWaarde}.

\section{TP representation of nonlinear systems}\label{SecPolySector}

To solve the controller synthesis problem from the previous section, we shortly recap the data-based polynomial representation of nonlinear functions based on TPs from \cite{MartinTP2}. According to Taylor's theorem \citep{TaylorRef}, we can write $f_i(x)=T_k[f_i(x)]+R_k[f_i(x)],i=1,\dots,n_x,$ with the TP of order $k$ at $x_{0}=0$
\begin{align*}
T_k[f_i(x)]=\sum_{|\alpha|=1}^{k}\frac{1}{\alpha !}\frac{\partial^{|\alpha|}f_i(0)}{\partial x^{\alpha}}x^\alpha={a_i^*}^Tz_i(x),
\end{align*}
where $\alpha!=\alpha_1!\cdots\alpha_{n_x}!$, the vector $z_i(x)\in\mathbb{R}[x]^{n_{z_i}}$ summarizes the polynomials $\frac{1}{\alpha !}x^\alpha$ and $a_i^*\in\mathbb{R}^{n_{z_i}}$ summarizes the unknown coefficients $\frac{\partial^{|\alpha|}f_i(0)}{\partial x^{\alpha}}$ for $|\alpha|=1,\dots,k$. Moreover, for all $x\in\mathbb{R}^{n_x}$ there exists a $\nu\in[0,1]$ such that $R_k[f_i(x)]$ is equal to the Lagrange remainder
\begin{align*}
\sum_{j=1}^{C_{k+n_x,n_x-1}}\frac{1}{\rho_j !}\frac{\partial^{k+1}f_i(\nu x)}{\partial x^{\rho_j}}x^{\rho_j},
\end{align*}
where $\rho_j,j=1,\dots,C_{k+n_x,n_x-1},$ correspond to all vectorial indices $\alpha$ with $|\alpha|=k+1$. While the results of this section also hold for TPs at arbitrary points $x_0\in\mathbb{R}^{n_x}$ as in \cite{MartinTP2}, the controller synthesis in Section~\ref{SecConSyn} is restricted to one TP at the equilibrium point $x=0$.

Since the existence of $\nu$ follows from the mean value theorem, its actual value is unknown. Therefore, \cite{MartinTP2} suggests two bounds on the remainder to circumvent the calculation of $\nu$ and the nonlinearity of the remainder.

\begin{lem}[\cite{MartinTP2}]\label{LemmaBounds}
Under Assumption~\ref{AssBoundDeri}, it holds true that $(R_k[f_i(x)])^2\leq R^{\text{abs}}_k[f_i(x)]\leq R^{\text{poly}}_k[f_i(x)]$ with
\begin{align}
&R^{\text{abs}}_k[f_i(x)]= \left(\sum_{j=1}^{C_{k+n_x,n_x-1}}\frac{M_{i,\rho_j}}{\rho_j !}\,||x^{\rho_j}||_2\right)^2,\label{AbsRemainderBound}\\
&R^{\text{poly}}_k[f_i(x)]=\kappa_i\sum_{j=1}^{C_{k+n_x,n_x-1}}\frac{M_{i,\rho_j}^2}{\rho_j !^2}x^{2\rho_j},\label{PolyRemainderBound}	
\end{align}
where $\kappa_i\in\mathbb{N}$ is the number of $M_{i,\rho_j}\neq0$ for $j=1,\dots,C_{k+n_x,n_x-1}$.
\end{lem}

Furthermore, due to Assumption~\ref{AssBoundDegree}, the $i$-th row of $B(x)$ can be written as ${b_i^*}^TG_i(x)$ with unknown coefficients $b_i^*$ and known polynomial matrix $G_i(x)$. Moreover, let $T_k[f(x)]=\begin{bmatrix}T_k[f_1(x)] & \cdots & T_k[f_{n_x}(x)]\end{bmatrix}^T$ and analogously for $R_k[f(x)]$ and let matrices $S_i$ suffice $\begin{bmatrix}{a_i^*}^T& {b_i^*}^T\end{bmatrix}^T=S_i\Theta^*$, i.e., $\Theta^*\in\mathbb{R}^{n_{\Theta}}$ summarizes all unknown coefficients of $T_k[f(x)]$ and $B(x)$. Finally, combining Taylor's theorem and \eqref{PolyRemainderBound} constitutes the polynomial description of \eqref{TrueSystem}
\begin{align}
f(x)+B(x)u=&T_k[f(x)]+R_k[f(x)]+B(x)u\notag\\
=&\begin{bmatrix}\begin{bmatrix}{a_1^*}^T & {b_1^*}^T\end{bmatrix}\begin{bmatrix}z_1(x)\\ G_1(x)u\end{bmatrix}\\\vdots\\ \begin{bmatrix}{a_{n_x}^*}^T & {b_{n_x}^*}^T\end{bmatrix}\begin{bmatrix}z_{n_x}(x)\\ G_{n_x}(x)u\end{bmatrix}\end{bmatrix}+R_k[f(x)]\label{StrucDyn2}\\
=&\underbrace{\begin{bmatrix}\begin{bmatrix}z_1(x)\\ G_1(x)u\end{bmatrix}^TS_1\\\vdots\\ \begin{bmatrix}z_{n_x}(x)\\ G_{n_x}(x)u\end{bmatrix}^TS_{n_x}\end{bmatrix}}_{=:Z(x,u)}\Theta^*+R_k[f(x)]\label{StrucDyn}
\end{align}
together with
\begin{align}\label{PolySecBound}
R_k[f(x)]^TR_k&[f(x)]\leq\omega(x)^T\mathcal{D}\omega(x),
\end{align}
where
\begin{align}
\mathcal{D}&=\sum_{i=1}^{n_x}\kappa_i\,\text{diag}\left(M_{i,\rho_1}^2\big| \cdots \big|M_{i,\rho_{C_{k+n_x,n_x-1}}}^2\right),\label{D_cal}\\
\omega(x)&=\begin{bmatrix}x^{\rho_1}/\rho_1!\\\vdots\\x^{\rho_{C_{k+n_x,n_x-1}}}/\rho_{C_{k+n_x,n_x-1}}!\end{bmatrix}.\label{omega}
\end{align}
Since the system description by \eqref{StrucDyn2} with \eqref{PolySecBound} is polynomial, a robust state-feedback design by SOS optimization is possible if system \eqref{TrueSystem} is known. Otherwise, a data-based inference on the unknown coefficients $\Theta^*$ is additionally required, which is shown in Section~\ref{SecBayesianSetMem}.

In contrast to \cite{Umenberger} and \cite{DePersis2} with $z_1=\dots=z_{n_x}$, we propose by \eqref{StrucDyn} a more flexible row-wise description with potentially distinct vectors $z_1,\dots,z_{n_x}$. This enables us to refine the accuracy of our data-based method by leveraging prior knowledge on the causality of the dynamics $f(x)$ and additional information on the polynomials of each element of $B(x)$. For instance, if $f_1$ is only a function of $x_1$ and the first row of $B(x)$ is constant, then this additional information can be incorporated by $z_1(x)=z_1(x_1)$ and $G_1(x)=G_1$. Moreover, if it is a priori known that $\begin{bmatrix}{a_i^*}^T& {b_i^*}^T\end{bmatrix}$ and $\begin{bmatrix}{a_j^*}^T& {b_j^*}^T\end{bmatrix}$ contain partially coinciding entries, then this redundancy can be considered by a reduced $\Theta^*$ and by the corresponding $S_i$ and $S_j$. For further examples, we refer to the numerical example in Section~\ref{SecNumEx} and Remark 2 in \cite{MartinTP1}. Note that the structure of \eqref{StrucDyn2} could also be incorporated in the general framework of \cite{Berberich} using a linear fraction representation with diagonal uncertainty description. However, since the unknown coefficients in \eqref{StrucDyn} are summarized in the vector $\Theta^*$, this row-wise consideration is preferable for their Bayesian inference in Section~\ref{SecBayes}.

\section{Gaussian inference of unknown coefficients $\Theta^*$}\label{SecBayesianSetMem}

In order to infer the unknown coefficients $\Theta^*$ in \eqref{StrucDyn} from data~\eqref{DataSystem}, we examine two approaches which are also known as Frequentist and Bayesian treatment \citep{Bishop} (Section 1.2).

\subsection{Frequentist perspective}\label{SecFrequ}

In the sequel, we show that a Frequentist treatment to infer $\Theta^*$ can be solved by the data-driven approaches with pointwise deterministic noise characterization. To this end, the following auxiliary result will be useful.

\begin{lem}[\cite{Cochran}]\label{LemChi}
Let $X\sim\mathcal{N}(\mu,\Xi)$ with $\mu\in\mathbb{R}^k$. Then, $Y=(X-\mu)^T\Xi^{-1}(X-\mu)$ is chi-squared distributed with $k$ degrees of freedom.	
\end{lem}

Since $\tilde{d}_i\sim\mathcal{N}(0,\sigma^2I),i=1,\dots,S,$ and are independent, Lemma~\ref{LemChi} implies
\begin{align*}
&p(\sigma^{-2}\tilde{d}_i^T\tilde{d}_i\leq Q_{n_x}(\delta^{\frac{1}{S}}),i=1,\dots,S)\\
=&\prod_{i=1}^{S}p(\sigma^{-2}\tilde{d}_i^T\tilde{d}_i\leq Q_{n_x}(\delta^{\frac{1}{S}}))=\delta.
\end{align*}
Therefore, the disturbance satisfies the noise description $\tilde{d}_i^T\tilde{d}_i\leq \sigma^2Q_{n_x}(\delta^{\frac{1}{S}}),i=1,\dots,S,$ with probability (w.p.) $\delta$. Together with $\tilde{d}_i=\dot{\tilde{x}}_i-Z(\tilde{x}_i,\tilde{u}_i)\Theta^*-R_k[f(\tilde{x}_i)]$, we derive the quadratic constraints $\star^TI\cdot(\dot{\tilde{x}}_i-Z(\tilde{x}_i,\tilde{u}_i)\Theta^*-R_k[f(\tilde{x}_i)])\leq \sigma^2Q_{n_x}(\delta^{\frac{1}{S}}),i=1,\dots,S,$ and proceed as in \cite{MartinTP2} (Section 3.B.) to conclude on a matrix $\varDelta_\text{F}\in\mathbb{R}^{(n_\Theta+1)\times(n_\Theta+1)}$ such that the set-membership 
\begin{align}\label{Sigma_Frequ}
\Sigma_{\text{F}}=\left\{\Theta:\begin{bmatrix}I\\\Theta^T\end{bmatrix}^T\varDelta_\text{F}\begin{bmatrix}I\\\Theta^T\end{bmatrix}\preceq0\right\}
\end{align}
contains $\Theta^*$ w.p. $\delta$. We refer to \cite{MartinTP1} for further insights and explain next that $\Sigma_{\text{F}}$ can be seen as a confidence region from a Frequentist treatment.

For that purpose, we compute $\dot{\tilde{x}}_i|\Theta\sim\mathcal{N}(Z(\tilde{x}_i,\tilde{u}_i)\Theta$ $ +R_k[f(\tilde{x}_i)],\sigma^2I)$. Hence, Lemma~\ref{LemChi} for these Gaussian random vectors results in the same conditions $\star^TI\cdot(\dot{\tilde{x}}_i-Z(\tilde{x}_i,\tilde{u}_i)\Theta^*-R_k[f(\tilde{x}_i)])\leq \sigma^2Q_{n_x}(\delta^{\frac{1}{S}}),i=1,\dots,S$. Concurrently, these conditions describe a confidence region of the conditional distribution
\begin{equation*}
p\left(\dot{\tilde{x}}_1,\dots,\dot{\tilde{x}}_S\big|\Theta\right)=\prod_{i=1}^{S}p\left(\dot{\tilde{x}}_i\big|\Theta\right)
\end{equation*}
which is typically considered in the Frequentist viewpoint \citep{Bishop}. Thus, if the generating process of data~\eqref{DataSystem} is repeated then the set-membership $\Sigma_{\text{F}}$ contains the true coefficients $\Theta^*$ in $100\cdot\delta$ percent of all repetitions. However, the problems arise that the set of data is only available once and $\Sigma_{\text{F}}$ might also be empty. Therefore, we also analyze the alternative Bayesian viewpoint. Note that we can also determine by
\begin{align*}
p\left(\sigma^{-2}\begin{bmatrix}\tilde{d}_1\\\vdots\\\tilde{d}_S\end{bmatrix}^T\begin{bmatrix}\tilde{d}_1\\\vdots\\\tilde{d}_S\end{bmatrix}\leq Q_{n_xS}(\delta)\right)=\delta
\end{align*}
a cumulatively bounded noise characterization, which resembles Section 6.C of \cite{DePersisNonlinearity}.

\subsection{Bayesian perspective}\label{SecBayes}

\begin{figure*}
\vspace{-0.2cm}
\setcounter{MYtempeqncnt}{\value{equation}}
\setcounter{equation}{12}
\begin{align}\label{EllipOuter}
-\begin{bmatrix}\Delta_{1} & 0 & \Delta_{2}\\ 0 & 0 & 0\\ \Delta_{2}^T & 0 & \Delta_3\end{bmatrix}+\begin{bmatrix}0 & 0 & 0\\ 0&  {W} & -{W}\mu_{\Theta1}\\ 0 & (-{W}\mu_{\Theta1})^T & \mu_{\Theta1}^T{W}\mu_{\Theta1}-\sum\nolimits_{i=1}^{n_\Theta}\eta_i\end{bmatrix}+\eta_0\begin{bmatrix}\Xi_\Theta^{-1} & -\Xi_\Theta^{-1} & 0\\ -\Xi_\Theta^{-1} & \Xi_\Theta^{-1} & 0\\ 0 & 0 & -Q_{n_\Theta}(\delta)\end{bmatrix}
\succeq0.
\end{align}
\setcounter{equation}{\value{MYtempeqncnt}}
\vspace{-0.2cm}
\hrulefill
\end{figure*}

While $\Theta^*$ is a deterministic vector in the Frequentist view, the true coefficients are a sample of a random vector $\Theta$ in the Bayesian treatment. To compute the distribution of $\Theta$, we update the a prior belief about the distribution by using the available data. For that reason, we first deduce a credibility region by adapting Proposition 2.1 from \cite{Umenberger}.

\begin{lem}\label{LemCredReg}
Let $\sum_{i=1}^{S}Z(\tilde{x}_i,\tilde{u}_i)^TZ(\tilde{x}_i,\tilde{u}_i)\succ0$ and the prior over the parameters $\Theta$ be uniform, i.e., $p(\Theta)\propto1$. Then, the posterior distribution is $p\left(\Theta\big|\dot{\tilde{x}}_1,\dots,\dot{\tilde{x}}_S\right){\sim}\mathcal{N}(\mu_{\Theta},\Xi_{\Theta})$ with $\mu_{\Theta}=\sigma^{-2}\Xi_{\Theta}\left(\sum_{i=1}^{S}Z(\tilde{x}_i,\tilde{u}_i)^T(\dot{\tilde{x}}_i-R_k[f(\tilde{x}_i)])\right)$ and $\Xi_{\Theta}^{-1}=\sigma^{-2}\sum_{i=1}^{S}Z(\tilde{x}_i,\tilde{u}_i)^TZ(\tilde{x}_i,\tilde{u}_i)$. Moreover, the true coefficients $\Theta^*$ are an element of the credibility region $\Sigma_{\text{Cred}}=\left\{\Theta:	\star^T\Xi_{\Theta}^{-1}\cdot\left(\Theta-\mu_{\Theta}\right)\leq Q_{n_{\Theta}}(\delta)\right\}$ w.p. $\delta$.	
\end{lem} 
\begin{pf}
Since \eqref{StrucDyn} is linear in the parameters $\Theta$ and the remainder is a deterministic vector, we can retrieve the posterior distribution by Bayes' rule as in \cite{Umenberger} (Proposition 2.1). The credibility region follows immediately by Lemma~\ref{LemChi} for the Gaussian posterior distribution.\hfill\quad \qed		
\end{pf}

Lemma~\ref{LemCredReg} supposes persistence of excitation of the data set \eqref{DataSystem}. It is not surprising that the calculation of the set-membership \eqref{Sigma_Frequ} by \cite{MartinTP2} (Proposition 1) requires the same assumption. Furthermore, the inclusion of $\Theta$ in $\Sigma_\text{Cred}$ in Lemma~\ref{LemCredReg} is given w.r.t. the posterior distribution $p\left(\Theta\big|\dot{\tilde{x}}_1,\dots,\dot{\tilde{x}}_S\right)$ whereas the inclusion in $\Sigma_\text{F}$ in the Frequentist treatment is regarding $p\left(\dot{\tilde{x}}_1,\dots,\dot{\tilde{x}}_S\big|\Theta\right)$. Thus, the probabilistic guarantee of both viewpoints can not be compared directly, and thereby both viewpoints are reasonable. We refer to \cite{Bishop} for a general discussion and to Section~\ref{SecNumEx} for a comparison in the context of data-driven controller synthesis.

The credibility region $\Sigma_{\text{Cred}}$ is not applicable so far because the mean $\mu_{\Theta}$ contains the unknown remainder $R_k[f(\tilde{x}_i)]$. Hence, we exploit the (tighter) bound on the remainder \eqref{AbsRemainderBound} in the following. While the explicit solution of \cite{Umenberger} (Lemma 3.1) is conceivable, it yields rather conservative results in particular due to the typically larger uncertainties of the coefficients of high order monomials. Thus, we propose an alternative next.

To this end, observe that the credibility region $\Sigma_{\text{Cred}}$ is intrinsically an ellipsoid with centre $\mu_{\Theta}=\mu_{\Theta1}+\mu_{\Theta2}$ for $\mu_{\Theta1}=\sigma^{-2}\Xi_{\Theta}\sum_{i=1}^{S}Z(\tilde{x}_i,\tilde{u}_i)^T\dot{\tilde{x}}_i$ and $\mu_{\Theta2}=-\sigma^{-2}\Xi_{\Theta}\sum_{i=1}^{S}Z(\tilde{x}_i,\tilde{u}_i)^TR_k[f(\tilde{x}_i)]$.
However, since only the bound \eqref{AbsRemainderBound} on the remainder is available, the mean $\mu_{\Theta}$ is another uncertainty set with known centre $\mu_{\Theta1}$ 
\begin{align}\label{Hyperrec}
\left\{\mu_{\Theta}=\mu_{\Theta1}+\mu_{\Theta2}: R_k[f(\tilde{x}_i)] \text{ with } \eqref{AbsRemainderBound}\right\}\subseteq\mu_{\Theta1}\oplus\mathcal{R}.
\end{align}
$\mathcal{R}$ is the hyperrectangle with centre zero, symmetric w.r.t. all axes, and edge lengths $\ell_i\geq0,i=1,\dots,n_{\Theta}$, with
\begin{equation*}
\begin{bmatrix}
\ell_1\\\vdots\\\ell_{n_{\Theta}}\end{bmatrix}=\frac{2}{\sigma^2}\sum_{i=1}^{S}|\Xi_\Theta Z(\tilde{x}_i,\tilde{u}_i)^T|_\text{ele}\begin{bmatrix}
\sqrt{R_k^{\text{abs}}[f_1(\tilde{x}_i)]}\\\vdots\\ \sqrt{R_k^{\text{abs}}[f_{n_x}(\tilde{x}_i)]}
\end{bmatrix}
\end{equation*} 
and $|\cdot|_\text{ele}$ taking the absolute value of each element of a matrix. The over approximation \eqref{Hyperrec} follows directly from the fact that $Mv\leq |M|_\text{ele}|v|_\text{ele}$ for any matrix $M\in\mathbb{R}^{n\times m}$ and vector $v\in\mathbb{R}^{m}$, where the inequality has to be understood elementwise. Combining \eqref{Hyperrec} and Lemma~\ref{LemCredReg}, we conclude that $\Theta^*$ is an element of
\begin{equation}\label{CredRegion2}
\tilde{\Sigma}_{\text{Cred}}=\left\{\Theta:
\Theta^T\Xi_{\Theta}^{-1}\Theta\leq Q_{n_{\Theta}}(\delta)\right\}\oplus\mu_{\Theta1}\oplus\mathcal{R}
\end{equation}
w.p. $\delta$, which does not require the evaluation of the remainder $R_k[f(x)]$. Note that the computation of lengths $\ell_i$ might be conservative as the sum over the approximation errors of all samples is considered. Hence, taking only the local data around $x=0$ into account might reduce the volume or the diameter of $\tilde{\Sigma}_{\text{Cred}}$. 

While $\tilde{\Sigma}_{\text{Cred}}$ is actually feasible for a controller synthesis by the full-block S-procedure \citep{SchererLMI}, we expect computationally demanding SOS optimization problems. Therefore, we suggest to first compute an ellipsoidal outer approximation of $\tilde{\Sigma}_{\text{Cred}}$.  

\begin{thm}\label{ThmElliOuter}
Let $\Sigma_{\text{Cred}}$ be given as in Lemma~\ref{LemCredReg} and $\ell_i>0,i=1,\dots,n_\Theta$. Then, there exist a positive definite matrix $\Delta_{1}\in\mathbb{R}^{n_\Theta\times n_\Theta}$, $\Delta_{2}\in\mathbb{R}^{n_\Theta}$, and scalars $\eta_0,\dots,\eta_{n_\Theta}\geq0$ solving \eqref{EllipOuter} \setcounter{equation}{13}
for $W=\text{diag}(\eta_{1}/(\ell_1/2)^2|\cdots|\eta_{n_\Theta}/(\ell_{n_\Theta}/2)^2)$ and $\Delta_3=\Delta_{2}^T\Delta_{1}^{-1}\Delta_{2}-1$. Moreover, for 
$\varDelta_\text{Cred}=\begin{bmatrix}-\Delta_{1} & \Delta_{2}\\ \Delta_{2}^T & -\Delta_3\end{bmatrix}^{-1}$, $\tilde{\Sigma}_{\text{Cred}}$ is a subset of 
\begin{align}\label{Sigma_A}
\bar{\Sigma}_{\text{Cred}}=\left\{\Theta:\begin{bmatrix}I\\\Theta^T\end{bmatrix}^T\varDelta_\text{Cred}\begin{bmatrix}I\\\Theta^T\end{bmatrix}\preceq0\right\}.
\end{align}
\end{thm}
\begin{pf}
At first, we show that \eqref{EllipOuter} has a solution. For that purpose, choose $\Delta_{1}=\frac{\eta_0}{2}\Xi_\Theta^{-1}$, $\Delta_{2}=0$, and $\eta_1=\dots=\eta_{n_\Theta}=\eta$. By choosing $\eta_0$ and $\eta$ small enough such that $\beta=1-(1-\mu_{\Theta1}^T\tilde{W}\mu_{\Theta1})n_\Theta\eta-Q_{n_\Theta}(\delta)\eta_0>0$ with $\tilde{W}=\frac{1}{n_\Theta\eta}W$, the first and the third diagonal block of \eqref{EllipOuter} are positive definite. Hence, the Schur complement can be applied twice to derive the equivalent condition $\tilde{W}-\frac{\eta_0}{n_\Theta\eta}\Xi_\Theta^{-1}-\frac{n_\Theta\eta}{\beta}\tilde{W}\mu_{\Theta1}(\tilde{W}\mu_{\Theta1})^T\succeq0$. Since $\frac{\eta_0}{\eta}\rightarrow0$ for $\eta\gg\eta_0$, $\frac{\eta}{\beta}\rightarrow0$ for $\eta,\eta_0\rightarrow0$, and $\tilde{W}$ is positive definite and independent of $\eta$ and $\eta_0$, we can always find sufficiently small $\eta,\eta_0>0$ satisfying \eqref{EllipOuter}. Second, we prove that $\tilde{\Sigma}_{\text{Cred}}\subseteq\bar{\Sigma}_{\text{Cred}}$. To this end, pre- and postmultiplying $\begin{bmatrix}\Theta^T & \mu_{\Theta}^T & 1\end{bmatrix}^T$ to \eqref{EllipOuter} implies 
\begin{align*}
&\begin{bmatrix}\Theta\\1\end{bmatrix}^T\begin{bmatrix}\Delta_{1} & \Delta_{2}\\ \Delta_{2}^T & \Delta_3\end{bmatrix}\begin{bmatrix}\Theta\\1\end{bmatrix}\\
-&\sum\nolimits_{i=1}^{n_\Theta}\eta_i\left((2/\ell_i)^2\left(\mu_\Theta[i]-\mu_{\Theta1}[i]\right)^2-1\right)\\
-&\eta_0\begin{bmatrix}\Theta-\mu_\Theta\\1\end{bmatrix}^T\begin{bmatrix}\Xi_\Theta^{-1} & 0\\ 0 & -Q_{n_\Theta}(\delta)\end{bmatrix}\begin{bmatrix}\Theta-\mu_\Theta\\1\end{bmatrix}\leq0,  
\end{align*}
where $\mu_\Theta[i]$ denotes the $i$-th element of $\mu_\Theta$ and respectively for $\mu_{\Theta1}[i]$. From the S-procedure related to \cite{Boyd} (page 46), we conclude that the ellipsoid
\begin{align}\label{DualElli}
\begin{bmatrix}\Theta\\1\end{bmatrix}^T\begin{bmatrix}\Delta_{1} & \Delta_{2}\\ \Delta_{2}^T & \Delta_3\end{bmatrix}\begin{bmatrix}\Theta\\1\end{bmatrix}\leq0  
\end{align}
comprises the Minkowski sum of the ellipsoid $\Theta^T\Xi_\Theta^{-1}\Theta\leq Q_{n_\Theta}(\delta)$ and the hyperrectangle $\mu_{\Theta1}\oplus\mathcal{R}$ which corresponds to \eqref{CredRegion2}. Furthermore, since $\Delta_{1}\succ0$ and the inverse
\begin{equation*}
\begin{bmatrix}\Delta_{1} & \Delta_{2}\\ \Delta_{2}^T & \Delta_3\end{bmatrix}^{-1}=\begin{bmatrix}\Delta_{1}^{-1}-\Delta_{1}^{-1}\Delta_{2}(\Delta_{1}^{-1}\Delta_{2})^T & \Delta_{1}^{-1}\Delta_{2}\\ (\Delta_{1}^{-1}\Delta_{2})^T & -1\end{bmatrix}
\end{equation*}
exists, the dualization lemma in Chapter 4.4.1 of \cite{SchererLMI} implies that \eqref{DualElli} is equivalent to \eqref{Sigma_A}.  \textcolor{white}{|}\hfill\quad \qed
\end{pf}

Note that the definiteness condition \eqref{EllipOuter} can be reformulated as an LMI by a Schur complement such that the computation of the ellipsoidal outer approximation $\bar{\Sigma}_{\text{Cred}}$ of $\tilde{\Sigma}_{\text{Cred}}$ with minimal volume or diameter is computationally tractable. In case $\ell_i=0$, choose $\ell_i=\epsilon_i>0$ to apply Theorem~\ref{ThmElliOuter}. Moreover, set-membership $\bar{\Sigma}_{\text{Cred}}$ is characterized as \eqref{Sigma_Frequ} in the Frequentist treatment and as for deterministic noise descriptions in \cite{MartinTP2}. Thus, the set-memberships $\bar{\Sigma}_{\text{Cred}}$ and $\Sigma_{\text{F}}$ together with the polynomial representation \eqref{StrucDyn2} and \eqref{PolySecBound} are admissible to verifying, among others, dissipativity properties of the unknown system \eqref{TrueSystem} using the framework of \cite{MartinTP2}.

\section{Data-driven state-feedback design for TP system representations}\label{SecConSyn}

First, we combine the data-based polynomial representation \eqref{StrucDyn2}, \eqref{PolySecBound}, and \eqref{Sigma_A}  together with the elaborated robust control framework of \cite{SchererLMI} to globally asymptotically stabilize the nonlinear system \eqref{TrueSystem}. Subsequently, we shortly discuss extensions, e.g., by desired performance criteria and by a local synthesis. These results also hold for the Frequentist treatment \eqref{Sigma_Frequ} and deterministic noise characterizations \citep{MartinTP2} due to the same set-membership characterization.

\subsection{Global stabilization}

\begin{figure*}
\vspace{-0.2cm}
\begin{align}\label{GlobSyn}
\Psi(x)=\star^T
\begin{bmatrix}\begin{array}{c|c|c}
\begin{matrix}0 & I\\ I & 0	\end{matrix} & 0 &0 \\\hline 
0 & \begin{matrix}-\tau_0I & 0\\ 0 & \tau_0\mathcal{D}^{-1}	\end{matrix} &0 \\\hline 
0 & 0 &\begin{matrix}
\mathcal{T}\otimes\Delta_{3} & -(\mathcal{T}\otimes\Delta_{2})^T\\-\mathcal{T}\otimes\Delta_{2} & \mathcal{T}\otimes\Delta_{1}
\end{matrix} 
\end{array}	\end{bmatrix}\cdot	
\begin{bmatrix}\hspace{0.05cm}\begin{matrix} 0 & -P\Omega^T & -\begin{bmatrix}
\begin{bmatrix}F_1P\\ G_1K\end{bmatrix}^TS_1 \cdots  \begin{bmatrix}F_{n_x}P\\ G_{n_x}K\end{bmatrix}^TS_{n_x}\end{bmatrix}\\I & 0 & 0  \\\hline
-\frac{\partial z}{\partial x}^T & 0 &0 \\ 0 & I & 0  \\\hline
-\frac{\partial z}{\partial x}^T & 0 & 0 \\ 0 & 0 & I 
\end{matrix} \hspace{0.05cm}\end{bmatrix}
\end{align}
\vspace{-0.2cm}
\begin{align}\label{GlobSynS}
\star^T
\begin{bmatrix}\begin{array}{c|c}
\begin{matrix}0 & I\\ I & 0	\end{matrix} & 0 \\\hline 
0 & \begin{matrix}-\tau_0I & 0\\ 0 & \tau_0\mathcal{D}^{-1}	\end{matrix}
\end{array}	\end{bmatrix}\cdot	
\begin{bmatrix}\hspace{0.05cm}\begin{matrix} 
-\begin{bmatrix}
\begin{bmatrix}F_1P\\ G_1K\end{bmatrix}^TS_1\Theta \cdots  \begin{bmatrix}F_{n_x}P\\ G_{n_x}K\end{bmatrix}^TS_{n_x}\Theta\end{bmatrix}\frac{\partial z}{\partial x}^T & -P\Omega^T \\ I & 0\\\hline
-\frac{\partial z}{\partial x}^T & 0\\ 0 & I
\end{matrix} \hspace{0.05cm}\end{bmatrix}\succ0
\end{align}
\vspace{-0.2cm}
\begin{align}\label{GlobSynPrimal}
\star^T
\begin{bmatrix}\begin{array}{c|c}
\begin{matrix}0 & I\\ I & 0	\end{matrix} & 0 \\\hline 
0 & \begin{matrix}-\frac{1}{\tau_0}I & 0\\ 0 & \frac{1}{\tau_0}\mathcal{D}	\end{matrix}
\end{array}	\end{bmatrix}\cdot	
\begin{bmatrix}\hspace{0.05cm}\begin{matrix} I & 0 \\\frac{\partial z}{\partial x}\begin{bmatrix}
\begin{bmatrix}F_1P\\ G_1K\end{bmatrix}^TS_1\Theta \cdots  \begin{bmatrix}F_{n_x}P\\ G_{n_x}K\end{bmatrix}^TS_{n_x}\Theta\end{bmatrix}^T & \frac{\partial z}{\partial x}\textcolor{white}{\bigg|}\\\hline
0 & I \\ \Omega P& 0 
\end{matrix} \hspace{0.05cm}\end{bmatrix}\prec0
\end{align}
\hrulefill
\vspace{-0.2cm}
\end{figure*}

By combining the LMI-based robust control framework of \cite{SchererLMI} with SOS relaxations, we achieve a convex optimization problem that yields a globally asymptotically stabilizing state-feedback law. For that reason, we introduce the Lyapunov function $z(x)^TP^{-1}z(x)$ with a vector of monomials $z(x)\in\mathbb{R}^{n_z}$  with $z(0)=0$ and $x=\begin{bmatrix}I & 0\end{bmatrix}z$. Hence, there exist matrices $F_i(x)$ with $z_i=F_{i}z,i=1,\dots,n_x,$ and a matrix $\Omega(x)$ such that $\omega=\Omega z$ from \eqref{omega}. Also note that without loss of generality $\mathcal{D}\succ0$ as if $\mathcal{D}$ from \eqref{D_cal} has a zero diagonal element, then the corresponding entry in $\mathcal{D}$ and $\omega$ can be omitted.

\begin{thm}\label{ThmGlobalSyn}
Let a non-empty set-membership $\bar{\Sigma}_{\text{Cred}}$ with $\varDelta_\text{Cred}=\begin{bmatrix}-\Delta_{1} & \Delta_{2}\\ \Delta_{2}^T & -\Delta_3\end{bmatrix}^{-1}$ as in Theorem~\ref{ThmElliOuter} be given. If there exist scalars $\epsilon,\epsilon_\tau>0$, a matrix $P\in\mathbb{R}^{n_z\times n_z}\succ0$, a polynomial matrix $K(x)\in\mathbb{R}[x]^{n_u\times n_z}$, and polynomials $\tau_0(x),\dots,\tau_{n_x}(x)$ such that $\tau_i-\epsilon_\tau\in\text{SOS}[x],i=0,\dots,n_x,$ and $\Psi(x)-\epsilon I\in\text{SOS}[x]^{n_z+C_{k+n_x,n_x-1}+n_{\Theta}n_x}$ for $\Psi(x)$ from \eqref{GlobSyn} with $\mathcal{T}=\text{diag}(\tau_1|\cdots|\tau_{n_x})$, then the equilibrium $x=0$ of \eqref{TrueSystem}
is globally asymptotically stable under the state feedback $u(x)=K(x)P^{-1}z(x)$ w.p. $\delta$.	
\end{thm}
\begin{pf}
To prove the statement, note that $\Psi(x)\succ0$ for all $x\in\mathbb{R}^{n_x}$ as $\Psi(x)-\epsilon I$ is SOS. Analogously, $\tau_i(x)>0,i=0,\dots,n_x$, which implies
\begin{equation*}
\star^T\begin{bmatrix}
\mathcal{T}\otimes\Delta_{3} & -(\mathcal{T}\otimes\Delta_{2})^T\\-\mathcal{T}\otimes\Delta_{2} & \mathcal{T}\otimes\Delta_{1}
\end{bmatrix}\cdot\begin{bmatrix}I\\-I\otimes\Theta\end{bmatrix}\preceq0
\end{equation*}
for all $\Theta\in\bar{\Sigma}_{\text{Cred}}$ due to \eqref{DualElli}. Therefore, by 
\begin{equation*}
\begin{bmatrix}
-\frac{\partial z}{\partial x}^T & 0 &0 \\ 0 & 0 & I
\end{bmatrix}\underbrace{\begin{bmatrix}I & 0\\ 0 & I\\ 
	(I\otimes\Theta)\frac{\partial z}{\partial x}^T & 0 
	\end{bmatrix}}_{=:\Gamma}=-\begin{bmatrix}
\begin{bmatrix}I \\-I\otimes\Theta\end{bmatrix}\frac{\partial z}{\partial x}^T & 0\end{bmatrix},
\end{equation*} 
by pre- and postmultiplying $\Psi(x)\succ0$ by $\Gamma$, and by the full-block S-procedure from \cite{SProcedure}, the condition \eqref{GlobSyn} implies \eqref{GlobSynS} for all $\Theta\in\bar{\Sigma}_{\text{Cred}}$.
Since $\tau_0(x)>0$, $1/\tau_0(x)$ exists and the dualization lemma \citep{SchererLMI} (Chapter 4.4.1) can be applied for \eqref{GlobSynS} which amounts to the equivalent condition \eqref{GlobSynPrimal}.
Pre- and postmultiplying the vector $\begin{bmatrix}z(x)^TP^{-1}\hspace{0.2cm}  R(x)^T\end{bmatrix}^T$ to \eqref{GlobSynPrimal} yield $0>z^TP^{-1}\frac{\partial z}{\partial x}\left\{Z(x,KP^{-1}z)\Theta+R(x)\right\}+\{\dots\}^T\frac{\partial z}{\partial x}^TP^{-1}z-\frac{1}{\tau_0}(R(x)^TR(x)-z^T\Omega^T\mathcal{D}\Omega z).$
Together with the S-procedure and the radially unbounded Lyapunov function $zP^{-1}z$, we conclude that the origin of all systems $\dot{x}=Z(x,u)\Theta+R(x)$ with $u(x)=K(x)P^{-1}z(x)$, $\Theta\in\bar{\Sigma}_{\text{Cred}}$, and $R(x)^TR(x)-z^T\Omega^T\mathcal{D}\Omega z\leq0$ are globally asymptotically stable. The theorem is proven as the unknown system \eqref{TrueSystem} is contained within this set of systems w.p. $\delta$. Indeed, the remainder $R_k[f(x)]$ suffices \eqref{PolySecBound} with $\omega=\Omega z$ and the coefficients of the TP $\Theta^*$ are an element of the set-membership $\bar{\Sigma}_{\text{Cred}}$ w.p. $\delta$.\hfill\quad \qed
\end{pf}

If $z_1(x)=\dots=z_{n_x}(x)$, then \eqref{GlobSyn} corresponds to a dual version of Theorem 2 of \cite{MartinTP2} for verifying dissipativity with supply rate $s(x,u)=0$ and for an unbounded state space. However, since the primal condition \eqref{GlobSynPrimal} is here not linear w.r.t. $P$ and $\tau_0$, we obtain by means of the dualization lemma the equivalent condition \eqref{GlobSynS} which is linear in the optimization variables. Thereby, condition \eqref{GlobSyn} can be solved by a computationally tractable SOS optimization. 

Furthermore, for $R_k[f(x)]=0$ and $z_1(x)=\dots=z_{n_x}(x)$, Theorem~\ref{ThmGlobalSyn} reduces to the special case of polynomial systems from Theorem 2 in \cite{DePersis1}. This result is also used in \cite{DePersis2} to locally asymptotically stabilize a nonlinear system by its TP approximation. By incorporating the remainder into Theorem~\ref{ThmGlobalSyn}, we can determine a polynomial state feedback that renders the origin globally asymptotically stable. Moreover, the LMI-based framework of \cite{SchererLMI} allows for even more possibilities as discussed in the next subsection.

\subsection{Extensions of Theorem~\ref{ThmGlobalSyn}}\label{SecExtensions}

This section presents possible extensions of Theorem~\ref{ThmGlobalSyn} by performance criteria, a local synthesis, the reduction of computational complexity, and leveraging prior knowledge on the value of certain coefficients.

Analogously to \cite{SchererLMI}, we introduce the performance input $w_p(t)$ and performance output $z_p(t)$
\begin{equation}\label{PerfSystem}
\begin{aligned}
\dot{x}&=f(x)+B(x)u+W_p(x)w_p,\\
z_p&=C(x)z(x)+D_u(x)u+D_w(x)w_p,
\end{aligned}
\end{equation}
with known polynomial matrices $W_p, C, D_u$, and $D_w$, and the notion of quadratic performance
\begin{equation}\label{QuadPerformance}
\int_{0}^{\infty}\star^TP_p\cdot\begin{bmatrix}w_p(t)\\z_p(t)\end{bmatrix}\text{d}t\leq\epsilon_p\int_{0}^{\infty}w_p(t)^Tw_p(t)\text{d}t,
\end{equation}
for some $\epsilon_p>0$ and performance matrix $P_p$ with $P_p^{-1}=\begin{bmatrix}\tilde{Q}_p & \tilde{S}_p\\ \tilde{S}_p^T & \tilde{R}_p
\end{bmatrix}$ and $\tilde{Q}_p\preceq0$. This comprises, among others, an $\mathcal{L}_2$-gain bound $\gamma>0$ on $w_p\mapsto z_p$ for $\tilde{Q}_p=-1/\gamma I, \tilde{S}_p=0,$ and $\tilde{R}_p=\gamma I$. Pursuing the arguments of Theorem~\ref{ThmGlobalSyn} and Chapter 8.1.2 of \cite{SchererLMI}, the performance channel $w_p\mapsto z_p$ of \eqref{PerfSystem} satisfies the performance \eqref{QuadPerformance} under the state feedback $u(x)=K(x)P^{-1}z(x)$ if a solution as in Theorem~\ref{ThmGlobalSyn} exists, but for
\begin{equation*}
\tilde{\Psi}(x)=\begin{bmatrix}\begin{array}{c|c}\Psi(x) & 0\\\hline 0 & 0	\end{array}	\end{bmatrix}+\begin{bmatrix}\begin{array}{c|c}\begin{matrix}
\frac{\partial z}{\partial x}W_p\tilde{Q}_{p}W_p^T\frac{\partial z}{\partial x}^T & 0 & 0\\ 0 & 0 & 0\\ 0 & 0 & 0
\end{matrix} & \begin{matrix}
\psi_{12}\\0\\0	\end{matrix}\\\hline \begin{matrix}\hspace{1.03cm}\psi_{12}^T\hspace{1.03cm}  &0&0\end{matrix}& \psi_{22}	\end{array}	\end{bmatrix}
\end{equation*} 
instead of $\Psi(x)$ with $\psi_{12}(x)=-(CP+D_uK)^T-\frac{\partial z}{\partial x}W_p(-\tilde{Q}_{p}D_w^T+\tilde{S}_{p})$ and $\psi_{22}(x)=D_w\tilde{Q}_{p}D_w^T-D_w\tilde{S}_{p}-\tilde{S}_{p}^TD_w+\tilde{R}_{p}$. Notice that if $\tilde{\Psi}(x)-\epsilon I$ is SOS then $\Psi(x)$ is SOS because $\tilde{Q}_p\preceq0$. Thus, the state feedback also constitutes a globally asymptotically stable closed loop by Theorem~\ref{ThmGlobalSyn}.     

While the presented controller synthesis achieves global closed-loop guarantees, Assumption~\ref{AssBoundDeri} can be conservative for $\mathbb{R}^{n_x}$ or only be valid for a compact and convex set $\mathbb{X}=\{x\in\mathbb{R}^{n_x}: w_i(x)\leq 0,\ w_i\in\mathbb{R}[x], i=1,\dots,n_w\}$ as in \cite{MartinTP2}. In this case, we can replace $\Psi$ in Theorem~\ref{ThmGlobalSyn} by
\begin{equation}\label{LocalSyn}
\tilde{\Psi}(x)+\sum_{i=1}^{n_w}\begin{bmatrix}\begin{array}{c|c}0 & 0\\\hline 0 & T_i(x)\end{array}	\end{bmatrix}w_i(x)
\end{equation}
for some to-be-optimized SOS matrices $T_i,i=1,\dots,n_w,$ to impose a state feedback that renders $x=0$ globally asymptotically stable whereas, the performance only holds for trajectories within $\mathbb{X}$. A related result can be found in \cite{SOSControllerDesign} but with full-block multipliers, $\tilde{\Psi}(x)+\sum_{i=1}^{n_w} T_i(x)w_i(x)$, that guarantees an asymptotically stable equilibrium, whose region of attraction contains the largest sublevel set of the Lyapunov function $z(x)^TP^{-1}z(x)$ within $\mathbb{X}$.


To reduce the computational complexity of Theorem~\ref{ThmGlobalSyn}, the Bayesian treatment of Section~\ref{SecBayes} can be employed to gather one set-membership for each $\begin{bmatrix}{a_i^*}^T & {b_i^*}^T\end{bmatrix}^T\in\mathbb{R}^{n_{\Theta_i}}, i=1,\dots,n_x$. By the S-procedure, these $n_x$ set-memberships can be considered in $\Psi(x)$ with $S_i=I$ and the multipliers $\tau_1,\dots,\tau_{n_x}$ as in Theorem~\ref{ThmGlobalSyn}. Thereby, the dimensions of $\Psi(x)$ reduce to $n_z+C_{k+n_x,n_x-1}+\sum_{i=1}^{n_x}n_{\Theta_i}$. 

In addition to the prior knowledge on the structure of the dynamics, we might have access to the actual value of certain coefficients. For instance, the $i$-th row corresponds to an integrator dynamics of a mechanical system or contains the TP of a known nonlinearity. Then, we could write the $i$-th row of \eqref{StrucDyn2} as
\begin{equation*}
\begin{bmatrix}a_{i,\text{prior}}^T & b_{i,\text{prior}}^T\end{bmatrix}\begin{bmatrix}z(x)\\ G_{i,\text{prior}}(x)u\end{bmatrix}+\begin{bmatrix}{a_i^*}^T & {b_i^*}^T\end{bmatrix}\begin{bmatrix}z_i(x)\\ G_i(x)u\end{bmatrix}
\end{equation*}
with everything known except for $\begin{bmatrix}{a_i^*}^T & {b_i^*}^T\end{bmatrix}$. This additive prior knowledge can be utilized in the procedure of Section~\ref{SecBayes} and in Theorem~\ref{ThmGlobalSyn} similar to \cite{Berberich}.

\section{Numerical Example}\label{SecNumEx}

This numerical example studies the stabilization of the unstable equilibrium $x=0$ of an inverted pendulum 
\begin{align}\label{InvertPen}
\dot{x}_1={a_1^*}^Tx,\quad \dot{x}_2=f_2(x_1)+b_2^*u,
\end{align}
with $a_1^*=\begin{bmatrix}0 & 1\end{bmatrix}^T$, $f_2(x_1)= \frac{g}{l}\sin(x_1)$, $b_2^*= \frac{1}{ml^2}$, $g=9.81, l=0.5$, and $m=0.2$. We assume that the structure in \eqref{InvertPen} is known but $a_1^*$, $f_2$, and $b_2^*$ are unknown. Furthermore, let the conservative upper bound $M_{2,[k+1\ 0]}=2\frac{g}{l},k\in\mathbb{N},$ be known which satisfies Assumption~\ref{AssBoundDeri} as $\Big|\Big|\frac{\partial^k f_2(x_1)}{\partial x_1^k}\Big|\Big|_2\leq\frac{g}{l},\forall x_1\in\mathbb{R}$. The inverted pendulum is numerically simulated for $10$ trajectories with random initial condition $x(0)\in[-0.1\ 0.1]^2$ and random but constant input signal $u(t)\in[-1\ 1]$. For $6$ samples from each trajectory with sampling time $0.1$, we evaluate the system dynamics \eqref{InvertPen} and add Gaussian noise with standard deviation $\sigma$, which constitute the data set \eqref{DataSystem} with a total of $S=60$ samples.

From the given data, we first calculate the set-memberships $\Sigma_{\text{F}}$ and $\bar{\Sigma}_{\text{Cred}}$ for the linear ($k=1$) and the third order ($k=3$) TP. Then, we apply Theorem~\ref{ThmGlobalSyn} with \eqref{LocalSyn} for $W_p=I$, $z_p=\begin{bmatrix}x^T&10^4u\end{bmatrix}^T$, $\tilde{Q}_p=-1/\gamma I, \tilde{S}_p=0,$ $\tilde{R}_p=\gamma I$, and the operating set $\mathbb{X}$ with $w_1(x)=x^Tx-1$. By minimizing the bound on the $\mathcal{L}_2$-gain $\gamma>0$ for a $z_p$ with such a large weighting of the control input, we achieve a state feedback with small control energy within $\mathbb{X}$. On the other hand, we could derive a closed loop that converges faster to the origin but requires more control energy by reducing the weighting of $u$. To globally stabilize the TP representation, we choose $z(x)=x$ and a feedback matrix $K(x)$ with degree $2$ and $4$ for $k=1$ and $k=3$, respectively.

\begin{figure}
\centering
\includegraphics[width=0.95\linewidth]{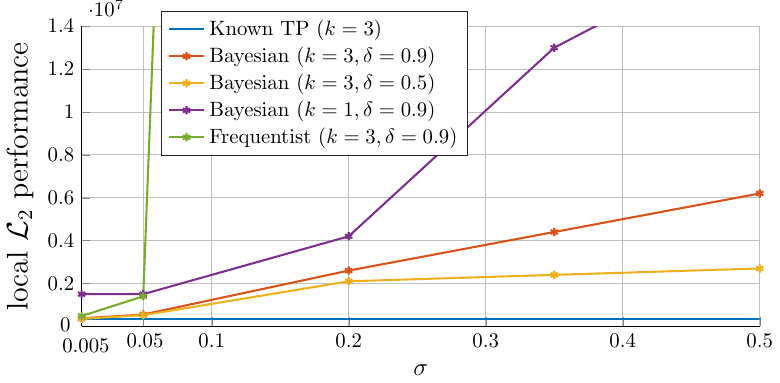}
\caption{Local performance of closed loops.}
\label{Fig.Performance}
\end{figure}

Figure~\ref{Fig.Performance} depicts the smallest bounds on the local $\mathcal{L}_2$ performance for credibility and confidence regions w.p. $\delta$. In comparison, a controller from Theorem~\ref{ThmGlobalSyn} with \eqref{GlobSyn} instead of \eqref{LocalSyn}, i.e., without optimized performance criterion, yields a local $\mathcal{L}_2$ performance of $2\cdot10^8$ for $k=3, \delta=0.9$, and $\sigma=0.2$.  
All attained closed loops exhibit global asymptotic stability despite the fact that we can only guarantee global asymptotic stability w.p. $\delta$. One explanation is that the derivation of the set-memberships encloses additional conservatism due to polynomial approximation error bounds. As expected, the performance of the data-driven controllers increases for larger $\sigma$ and is comparable for small $\sigma$ to the performance of a controller which is derived from Theorem~\ref{ThmGlobalSyn} with \eqref{LocalSyn} and system representation \eqref{StrucDyn} with known third order TP and \eqref{PolySecBound}. We assess the Frequentist treatment to be excessively conservative w.r.t. to Gaussian noise compared to the Bayesian treatment. The optimization problem from Theorem~\ref{ThmGlobalSyn} with \eqref{LocalSyn} is solved for all scenarios by YALMIP with solver MOSEK in Matlab in less than $8\,\text{s}$ on a Lenovo i5 notebook.

\section{Conclusion}

Within the framework of the data-based Taylor polynomial representation of general nonlinear systems from \cite{MartinTP2}, we investigated 
a Frequentist and a Bayesian treatment for Gaussian inference of the underlying unknown Taylor polynomial. Moreover, we combined this result with the robust control framework \citep{SchererLMI} for a data-driven controller synthesis by SOS optimization to determine state-feedback laws that render a known equilibrium globally asymptotically stable while satisfying a (local) quadratic performance. Our results can be combined with prior knowledge on the dynamics to improve accuracy and data-efficiency.

Contrary to the presented indirect controller synthesis, interesting future work includes a direct controller design by the full-block S-procedure and \cite{SchererSOS}. Furthermore, to reduce the conservatism of the controller due to the polynomial approximation, one interesting extension of Theorem~\ref{ThmGlobalSyn} is the consideration of multiple approximation polynomials as a piecewise polynomial representation, as suggested by \cite{MartinTP1}.


\end{document}